\documentclass[a4paper,12pt,twoside,leqno,final]{amsart}
\usepackage{amsmath}
\usepackage{amssymb}

\newtheorem{thm}{Theorem}[section]
\newtheorem{lem}[thm]{Lemma}

\newtheorem{prop}[thm]{Proposition}

\newtheorem*{tha}{Theorem A}
\newtheorem*{thb}{Theorem B}

\newcommand{\C}{{\mathbb C}}
\newcommand{\D}{{\mathbb D}}
\newcommand{\R}{{\mathbb R}}
\newcommand{\T}{{\mathbb T}}
\newcommand{\Z}{{\mathbb Z}}
\newcommand{\N}{{\mathbb N}}

\newcommand{\bmo}{{\rm BMO}}
\newcommand{\bmoa}{{\rm BMOA}}

\newcommand{\La}{\Lambda}

\newcommand{\eps}{\varepsilon}

\newcommand{\wt}{\widetilde}
\newcommand{\f}{\frac}
\newcommand{\ov}{\overline}
\newcommand{\al}{\alpha}

\newcommand{\ze}{\zeta}
\renewcommand{\th}{\theta}

\newcommand{\ph}{\varphi}

\newcommand{\Om}{\Omega}
\newcommand{\Omte}{\Om(\th,\eps)}

\numberwithin{equation}{section}

\title[Interpolation and duality]
{Interpolation and duality in spaces\\ 
of pseudocontinuable functions}
\author{Konstantin M. Dyakonov}
\address{Departament de Matem\`atiques i Inform\`atica, IMUB, BGSMath, Universitat de Barcelona, Gran Via 585, E-08007 Barcelona, Spain}
\address{ICREA, Pg. Llu\'is Companys 23, E-08010 Barcelona, Spain}
\email{konstantin.dyakonov@icrea.cat}
\keywords{Hardy space, smoothness class, BMO, inner function, interpolating Blaschke product, star-invariant subspace, duality}
\subjclass[2010]{30H05, 30H10, 30J05, 46J15} 
\thanks{Supported in part by grants MTM2017-83499-P and PID2021-123405NB-I00 from El Ministerio de Ciencia e Innovaci\'on (Spain).}

\begin{document}
\begin{abstract}
Given an inner function $\theta$ on the unit disk, let $K^p_\theta:=H^p\cap\theta\overline z\overline{H^p}$ be the associated star-invariant subspace of the Hardy space $H^p$. Also, we put $K_{*\theta}:=K^2_\theta\cap{\rm BMO}$. Assuming that $B=B_{\mathcal Z}$ is an interpolating Blaschke product with zeros $\mathcal Z=\{z_j\}$, we characterize, for a number of smoothness classes $X$, the sequences of values $\mathcal W=\{w_j\}$ such that the interpolation problem $f\big|_{\mathcal Z}=\mathcal W$ has a solution $f$ in $K^2_B\cap X$. Turning to the case of a general inner function $\theta$, we further establish a non-duality relation between $K^1_\theta$ and $K_{*\theta}$. Namely, we prove that the latter space is properly contained in the dual of the former, unless $\theta$ is a finite Blaschke product. From this we derive an amusing non-interpolation result for functions in $K_{*B}$, with $B=B_{\mathcal Z}$ as above.
\end{abstract}

\maketitle

\section{Introduction and results}

We write $\T$ for the unit circle $\{\ze\in\C:|\ze|=1\}$ and $m$ for the normalized arc length measure on $\T$; thus, $dm(\ze)=|d\ze|/(2\pi)$. We then define the spaces $L^p:=L^p(\T,m)$ in the usual way and let $\|\cdot\|_p$ denote the standard norm on $L^p$. Also, for $1\le p\le\infty$, we introduce the {\it Hardy space} $H^p$ by putting 
$$H^p:=\{f\in L^p:\,\widehat f(n)=0\,\,\,\text{\rm for }\,n=-1,-2,\dots\},$$
where $\widehat f(n)$ is the $n$th {\it Fourier coefficient} of $f$ given by 
$$\widehat f(n):=\int_\T\ov\ze^nf(\ze)\,dm(\ze),\qquad n\in\Z.$$
The Poisson integral (i.e., harmonic extension) of an $H^p$ function being holomorphic on the disk
$$\D:=\{z\in\C:\,|z|<1\}$$
(see \cite[Chapter II]{G}), we may use this extension to view elements of $H^p$ as holomorphic functions on $\D$ when convenient. 
\par Furthermore, we write $P_+$ (resp., $P_-$) for the orthogonal projection from $L^2$ onto $H^2$ (resp., onto $\ov z\ov{H^2}=L^2\ominus H^2$). By a classical theorem of M. Riesz (see \cite[Chapter III]{G}), each of these projections admits a bounded extension---or restriction---to $L^p$, with $1<p<\infty$, and maps $L^p$ onto $H^p$ (resp., onto $\ov z\ov{H^p}$).

\par Now suppose $\th$ is an {\it inner function}, meaning that $\th\in H^\infty$ and $|\th|=1$ a.e. on $\T$. The corresponding {\it star-invariant} (or {\it model}) {\it subspace} $K^p_\th$ is then defined by 
\begin{equation}\label{eqn:defnkptheta}
K^p_\th:=\{f\in H^p:\,\ov z\ov f\th\in H^p\},\qquad1\le p\le\infty,
\end{equation}
so that $K^p_\th=H^p\cap\th\ov z\ov{H^p}$. (When $p=2$, yet another equivalent definition is $K^2_\th=H^2\ominus\th H^2$.) It is clear from \eqref{eqn:defnkptheta} that the antilinear isometry 
\begin{equation}\label{eqn:antilinisom}
f\mapsto\ov z\ov f\th=:\widetilde f
\end{equation}
leaves $K^p_\th$ invariant. Also, it is well known (see \cite{DSS, N}) that each $K^p_\th$ is invariant under the {\it backward shift operator}
$$\mathfrak B:f\mapsto\f{f-f(0)}z,\qquad f\in H^p,$$ 
and conversely, that every closed and nontrivial $\mathfrak B$-invariant subspace of $H^p$, with $1\le p<\infty$, arises in this way. 
\par The functions belonging to some $K^p_\th$ space (i.e., noncyclic vectors of $\mathfrak B$) are known as {\it pseudocontinuable} functions. In fact, they are characterized by the property of having a meromorphic pseudocontinuation to $\D_-:=\C\setminus(\D\cup\T)$; that is, the function in question should agree a.e. on $\T$ with the boundary values of some meromorphic function of bounded characteristic in $\D_-$ (see \cite{DSS} for details). 

\par The orthogonal projection from $H^2$ onto $K^2_\th$ is given by $f\mapsto\th P_-(\ov\th f)$, and the M. Riesz theorem shows that the same formula provides, for $1<p<\infty$, a bounded projection from $H^p$ onto $K^p_\th$ parallel to $\th H^p$. This yields the direct sum decomposition 
\begin{equation}\label{eqn:dirsumtheta}
H^p=K^p_\th\oplus\th H^p,\qquad1<p<\infty,
\end{equation}
with orthogonality for $p=2$. 

\par Among the inner functions $\th$, of special relevance to us are Blaschke products. Recall that, for a sequence $\mathcal Z=\{z_j\}\subset\D$ with
\begin{equation}\label{eqn:blacond}
\sum_j(1-|z_j|)<\infty, 
\end{equation}
the associated {\it Blaschke product} is given by 
$$B(z)=B_{\mathcal Z}(z):=\prod_j\f{|z_j|}{z_j}\f{z_j-z}{1-\ov z_jz}$$
(if $z_j=0$, then we set $|z_j|/z_j=-1$). The product converges uniformly on compact subsets of $\D$ and defines an inner function that vanishes precisely at the $z_j$'s; see \cite[Chapter II]{G}. If, in addition, 
\begin{equation}\label{eqn:intseq}
\inf_j|B'(z_j)|\,(1-|z_j|)>0,
\end{equation}
then we say that $B$ is an {\it interpolating Blaschke product}. Accordingly, the sequences $\mathcal Z=\{z_j\}$ in $\D$ that satisfy \eqref{eqn:blacond} and \eqref{eqn:intseq}, with $B=B_{\mathcal Z}$, are called {\it interpolating} (or $H^\infty$-{\it interpolating}) {\it sequences}. By a celebrated theorem of Carleson (see \cite{Carl} or \cite[Chapter VII]{G}), these are precisely the sequences $\mathcal Z$ with the property that 
$$H^\infty\big|_{\mathcal Z}=\ell^\infty.$$

\par Here and below, the following standard notation (and terminology) is used. Given a sequence $\mathcal Z=\{z_j\}$ of pairwise distinct points in $\D$, the {\it trace} $f\big|_{\mathcal Z}$ of a function $f:\D\to\C$ is defined to be the sequence $\{f(z_j)\}$; and if $\mathcal X$ is a certain function space on $\D$, then the corresponding {\it trace space} is 
$$\mathcal X\big|_{\mathcal Z}:=\left\{f\big|_{\mathcal Z}:\,f\in\mathcal X\right\}.$$

\par We shall be concerned with interpolation problems for functions in star-invariant subspaces---specifically, for those in $K^p_B$, where $B$ is an interpolating Blaschke product. 
Some of the earlier results in this area can be found in \cite{AH, DPAMS, GW, HNP}, while others, more relevant to our current topic, will be recalled presently. 

\par First, we need yet another piece of notation. Given numbers $p>0$, $\gamma\in\R$ and a sequence $\mathcal Z=\{z_j\}\subset\D$, we write $\ell^p_\gamma(\mathcal Z)$ for the set of all sequences $\{w_j\}\subset\C$ satisfying
$$\sum_j|w_j|^p(1-|z_j|)^{\gamma}<\infty.$$
Now, if $1<p<\infty$ and if $B=B_{\mathcal Z}$ is an interpolating Blaschke product with zero sequence $\mathcal Z$, then we have 
$$K^p_B\big|_{\mathcal Z}=H^p\big|_{\mathcal Z}=\ell^p_1(\mathcal Z).$$
Indeed, the left-hand equality follows from \eqref{eqn:dirsumtheta} with $\th=B$, while the other holds by a well-known theorem of Shapiro and Shields \cite{SS}. In addition, for each sequence $\mathcal W=\{w_j\}$ in $\ell^p_1(\mathcal Z)$, there is a {\it unique} function $f\in K^p_B$ with $f\big|_{\mathcal Z}=\mathcal W$; the uniqueness is due to the fact that $K^p_B\cap BH^p=\{0\}$.

\par The case of $K^\infty_B$ is subtler, as the next result shows.

\begin{tha} Suppose that $\mathcal Z=\{z_j\}$ is an interpolating sequence in $\D$ and $B=B_{\mathcal Z}$ is the associated Blaschke product. Then we have
\begin{equation}\label{eqn:kinfbzlinf}
K^\infty_B|_{\mathcal Z}=\ell^\infty
\end{equation}
if and only if
\begin{equation}\label{eqn:unifc}
\sup\left\{\sum_j\f{1-|z_j|}{|\ze-z_j|}:\,\ze\in\T\right\}<\infty.
\end{equation}
\end{tha}

This theorem is essentially a consequence of Hru\v s\v cev and Vinogradov's work in \cite{HV}; see also \cite[Section 3]{C2} for details.

\par Condition \eqref{eqn:unifc} above is known as the {\it uniform Frostman condition}, and the sequences $\mathcal Z=\{z_j\}$ in $\D$ that obey it are called {\it Frostman sequences}. While a Frostman sequence need not be interpolating (in fact, its points are not even supposed to be pairwise distinct), it does necessarily split into finitely many interpolating sequences; see \cite{HV} for a proof. Finally, a Blaschke product whose zeros form a Frostman sequence will be referred to as a {\it Frostman Blaschke product}.

\par We mention in passing that, by a theorem of Vinogradov \cite{V}, the identity 
$$K^\infty_{B^2}|_{\mathcal Z}=\ell^\infty$$
is valid whenever $\mathcal Z$ is an interpolating sequence and $B=B_{\mathcal Z}$. It should be noted, however, that $K^\infty_{B^2}$ is strictly larger than $K^\infty_B$. 

\par To describe the trace class $K^\infty_B|_{\mathcal Z}$ in the general case (i.e., when \eqref{eqn:unifc} no longer holds), we first introduce a bit of notation. Once the interpolating sequence $\mathcal Z=\{z_j\}$ is fixed, we associate with each sequence $\mathcal W=\{w_j\}$ from $\ell^1_1(\mathcal Z)$ the {\it conjugate sequence} $\wt{\mathcal W}=\{\wt w_k\}$ whose elements are
\begin{equation}\label{eqn:conjseq}
\widetilde w_k:=\sum_j\f{w_j}{B'(z_j)\cdot(1-z_j\ov z_k)}\qquad(k=1,2,\dots).
\end{equation}
The absolute convergence of the series in \eqref{eqn:conjseq} is ensured, for any $k\in\N$, by the fact that $\mathcal W\in\ell^1_1(\mathcal Z)$ in conjunction with \eqref{eqn:intseq}. Because $\ell^1_1(\mathcal Z)$ contains $\ell^\infty$, as well as every $\ell^p_1(\mathcal Z)$ with $1<p<\infty$, the sequence $\wt{\mathcal W}$ is well defined whenever $\mathcal W$ belongs to one of these spaces.

\par The following result was established in \cite{DBLMS18}. 

\begin{thb} Suppose that $\mathcal Z=\{z_j\}$ is an interpolating sequence in $\D$ and $B=B_{\mathcal Z}$ is the associated Blaschke product. Given a sequence $\mathcal W\in\ell^\infty$, one has $\mathcal W\in K^\infty_B\big|_{\mathcal Z}$ if and only if $\wt{\mathcal W}\in\ell^\infty$.
\end{thb}

It was further conjectured in \cite{DBLMS18, DIEOT} that the trace space $K^1_B\big|_{\mathcal Z}$ is describable in similar terms, i.e., that the necessary conditions $\mathcal W\in\ell^1_1(\mathcal Z)$ and $\wt{\mathcal W}\in\ell^1_1(\mathcal Z)$ are also sufficient for $\mathcal W$ to be in $K^1_B\big|_{\mathcal Z}$. To the best of our knowledge, the conjecture is still open. 

\par Here, our purpose is to supplement Theorem B by characterizing the values of {\it smooth}, not just bounded, functions in $K^2_B$ on the (interpolating) sequence $\mathcal Z=B^{-1}(0)$. To be more precise, of concern are trace spaces of the form $\left(K^2_B\cap X\right)\big|_{\mathcal Z}$, where $X$ is a certain smoothness class on $\T$. Specifically, $X$ will be one of the following spaces. 

\smallskip$\bullet$ The {\it Lipschitz--Zygmund space} $\La^\al=\La^\al(\T)$ with $\al>0$. This is the set of functions $f\in C(\T)$ satisfying
$$\|\Delta_h^nf\|_\infty=O(|h|^\al),\qquad h\in\R,$$
where $n$ is some (any) integer with $n>\al$, and $\Delta_h^n$ denotes the $n$th order difference operator with step $h$. (As usual, the difference operators $\Delta_h^k$ are defined inductively: we put $$(\Delta_h^1f)(\ze):=f(e^{ih}\ze)-f(\ze),\qquad\ze\in\T,$$
and $\Delta_h^kf:=\Delta_h^1\Delta_h^{k-1}f$ for $k\ge2$.)

\smallskip$\bullet$ $\bmo=\bmo(\T)$, the space of functions of {\it bounded mean oscillation} on $\T$. Recall that an integrable function $f$ on $\T$ belongs to $\bmo$ if and only if 
$$\|f\|_*:=\left|\int_\T f\,dm\right|+\sup_I\f1{m(I)}\int_I|f-f_I|\,dm<\infty,$$
where $f_I:=m(I)^{-1}\int_I f\,dm$; the supremum is taken over the open arcs $I\subset\T$. Even though $\bmo$ contains discontinuous and unbounded functions, there are reasons for viewing it as a smoothness class. In a sense, it corresponds to the endpoint as $\al\to0$ of the $\La^\al$ scale. We also need the analytic subspace $\bmoa:=\bmo\cap H^2$. 

\smallskip$\bullet$ The {\it Gevrey class} $G_\al=G_\al(\T)$ with $\al>0$. This is the set of 
functions $f\in C^\infty(\T)$ satisfying 
$$\|f^{(n)}\|_\infty\le Q_f^{n+1}(n!)^{1+1/\al},\qquad n=0,1,2,\dots,$$
with some constant $Q_f>0$. Here, we write $f^{(n)}(e^{it})$ for the $n$th order derivative of the function $t\mapsto f(e^{it})$, which is assumed to be $C^\infty$-smooth on $\R$. 

\smallskip$\bullet$ The {\it Sobolev space} $\mathcal L^p_s=\mathcal L^p_s(\T)$ with $1<p<\infty$ and $s>0$, defined by
$$\mathcal L^p_s:=\{f\in L^p:\,f^{(s)}\in L^p\},$$
with the appropriate interpretation of the (possibly fractional) derivative $f^{(s)}$. Precisely speaking, we write $f^{(s)}\in L^p$ to mean that there is a function $g\in L^p$ satisfying $\widehat g(n)=(in)^s\widehat f(n)$ for all $n\in\Z$. 

\smallskip For each of these choices of $X$, we now characterize the sequences $\mathcal W$ from the trace space $\left(K^2_B\cap X\right)\big|_{\mathcal Z}$ in terms of the conjugate sequence $\wt{\mathcal W}$, as defined by \eqref{eqn:conjseq} above. The description always involves a certain decay condition (or growth restriction) on $\wt{\mathcal W}$, as we shall presently see.

\begin{thm}\label{thm:intsmooth} Let $\al>0$, $1<p<\infty$ and $s>0$. Also, let $X$ be one of the following spaces: $\La^\al$, $\bmo$, $G_\al$ or $\mathcal L^p_s$. Given an interpolating Blaschke product $B=B_{\mathcal Z}$ with zeros $\mathcal Z=\{z_k\}$ and a sequence $\mathcal W=\{w_k\}\in\ell^2_1(\mathcal Z)$, we have 
$$\mathcal W\in\left(K^2_B\cap X\right)\big|_{\mathcal Z}$$
if and only if 
\par{\rm (a)} $|\wt w_k|=O\left((1-|z_k|)^\al\right)$ when $X=\La^\al$, 
\par{\rm (b)} $\wt{\mathcal W}\in\ell^\infty$ when $X=\bmo$, 
\par{\rm (c)} there is a constant $c>0$ such that 
$$|\widetilde w_k|=O\left(\exp\left(-\f c{(1-|z_k|)^\al}\right)\right)$$ 
when $X=G_\al$,
\par{\rm (d)} $\wt{\mathcal W}\in\ell^p_{1-sp}(\mathcal Z)$ when $X=\mathcal L^p_s$.
\end{thm}

The intersection $K^2_B\cap\bmo$, which corresponds to case (b) above, will be henceforth denoted by $K_{*B}$. Similarly, for a general inner function $\th$, we define
$$K_{*\th}:=K^2_\th\cap\bmo.$$
Comparing Theorem B with the $\bmo$ part of Theorem \ref{thm:intsmooth}, we see that the structure of the trace space $K^\infty_B\big|_{\mathcal Z}$ is remarkably similar to that of $K_{*B}\big|_{\mathcal Z}$. In light of this observation, we may wonder what the $\bmo$ counterpart of Theorem A could look like. Specifically, we may ask if there exist infinite Blaschke products $B=B_{\mathcal Z}$ for which the trace space $K_{*B}\big|_{\mathcal Z}$ is completely determined by the natural (and necessary) logarithmic growth condition on the values. 

\par To be more precise, suppose that $\mathcal Z=\{z_k\}$ is a sequence of pairwise distinct points in $\D$, and write $\ell^\infty_{\log}(\mathcal Z)$ for the space of sequences $\mathcal W=\{w_k\}\subset\C$ with
$$|w_k|=O\left(\log\f2{1-|z_k|}\right).$$
It is well known (and easily shown) that every $f\in\bmoa$ satisfies 
$$|f(z)|=O\left(\log\f2{1-|z|}\right),\qquad z\in\D,$$
so $\bmoa\big|_{\mathcal Z}$ is always contained in $\ell^\infty_{\log}(\mathcal Z)$. The equality 
\begin{equation}\label{eqn:bmoazlog}
\bmoa\big|_{\mathcal Z}=\ell^\infty_{\log}(\mathcal Z)
\end{equation}
obviously need not hold in general, but it does actually occur for some infinite sequences $\mathcal Z=\{z_k\}$ (which form a tiny subfamily among the $H^\infty$-interpolating sequences). For instance, \eqref{eqn:bmoazlog} will be valid provided that 
$$|z_j-z_k|\ge c(1-|z_j|)^s,\qquad j\ne k,$$
for some constants $c>0$ and $s\in(0,\f12)$; see \cite[Theorem 11]{DCAG}. 

\par The question is what happens to \eqref{eqn:bmoazlog} when $\bmoa$ gets replaced by its subspace $K_{*B}$, with $B=B_{\mathcal Z}$. The property that arises is thus
\begin{equation}\label{eqn:kstarbzlog}
K_{*B}\big|_{\mathcal Z}=\ell^\infty_{\log}(\mathcal Z),
\end{equation}
and we regard it as an analogue of \eqref{eqn:kinfbzlinf} in the $\bmo$ setting. In contrast to \eqref{eqn:kinfbzlinf}, however, \eqref{eqn:kstarbzlog} does not lead to any nontrivial class of sequences. Indeed, our next result shows that \eqref{eqn:kstarbzlog} is only possible when $\mathcal Z=B^{-1}(0)$ is a finite set. 

\begin{prop}\label{prop:nonint} Whenever $B=B_{\mathcal Z}$ is an infinite Blaschke product with simple zeros, the trace space $K_{*B}\big|_{\mathcal Z}$ is properly contained in $\ell^\infty_{\log}(\mathcal Z)$.
\end{prop}

This will be deduced from another result, which deals with the case of a general inner function $\th$ and asserts an amusing lack of duality between the star-invariant subspaces $K^1_\th$ and $K_{*\th}$. 

\par It is well known that, for $1<p<\infty$, the dual of the Hardy space $H^p$ (under the pairing $\langle f,g\rangle=\int_\T f\ov g\,dm$) is $H^q$ with $q=p/(p-1)$, while the dual of $H^1$ is $\bmoa$; see, e.g., \cite[Chapter VI]{G}. The former duality relation has a natural counterpart in the $K^p_\th$ setting, namely $(K^p_\th)^*=K^q_\th$ for $p$ and $q$ as above (see \cite[Lemma 4.2]{C1}), and one may wonder if the identity $(K^1_\th)^*=K_{*\th}$ has any chance of being true, at least for some inner functions $\th$. Our last theorem says that this is never the case, except when $\th$ is a finite Blaschke product.

\begin{thm}\label{thm:nondual} Given an inner function $\th$, other than a finite Blaschke product, there exists a non-$\bmo$ function $g\in K^2_\th$ such that the functional
$$f\mapsto\int_\T f\ov g\,dm,$$
defined initially for $f\in K^2_\th$, extends continuously to $K^1_\th$.
\end{thm}

In other words, whenever $\th$ is an \lq\lq interesting" (i.e., nonrational) inner function, $K_{*\th}$ is properly contained in $(K^1_\th)^*$. One might compare this non-duality result with 
Bessonov's duality theorem for $K^1_\th$ that appears in \cite{B}. There, $\th$ was assumed to be a one-component inner function, meaning that the set $\{z\in\D:|\th(z)|<\eps\}$ is connected for some $\eps\in(0,1)$, and the dual of $K^1_\th\cap zH^1$ was identified with a certain discrete $\bmo$ space on $\T$. 

\par In the remaining part of the paper, we first list a number of auxiliary facts (these are collected in Section 2) and then use them to prove our current results. The proofs are in Sections 3 and 4. 

\smallskip{\it Acknowledgement.} I thank Carlo Bellavita for a helpful conversation. In particular, Theorem \ref{thm:nondual} of this paper arose in response to a question he asked me.

\section{Preliminaries}

Several background results will be needed. When stating the first of these, we shall assume that $X$ is one of our smoothness spaces (namely, $\La^\al$, $\bmo$, $G_\al$ or $\mathcal L^p_s$), the admissible values of the parameters $\al$, $p$ and $s$ being as above. 

\begin{lem}\label{lem:firstlem} Let $f\in H^2$ and let $B=B_{\mathcal Z}$ be an interpolating Blaschke product with zeros $\mathcal Z=\{z_k\}$. In order that $P_-(\ov Bf)\in X$, it is necessary and sufficient that 
\par{\rm (a)} $|f(z_k)|=O\left((1-|z_k|)^\al\right)$ when $X=\La^\al$, 
\par{\rm (b)} $\{f(z_k)\}\in\ell^\infty$ when $X=\bmo$, 
\par{\rm (c)} for some $c>0$, 
$$|f(z_k)|=O\left(\exp\left(-\f c{(1-|z_k|)^\al}\right)\right)$$ 
when $X=G_\al$,
\par{\rm (d)} $\{f(z_k)\}\in\ell^p_{1-sp}(\mathcal Z)$ when $X=\mathcal L^p_s$.
\end{lem}

The statements corresponding to parts (a) and (b) were proved in \cite{DSpb93} as Theorems 4.1 and 5.2. For parts (c) and (d), we refer to \cite{DIUMJ94}; specifically, see Theorems 1 and 7 therein.

\par Another (well-known) fact to be used below is that the space $\bmoa$ enjoys the so-called {\it K-property} of Havin, as defined in \cite{H}. The precise meaning of this assertion is as follows.

\begin{lem}\label{lem:sesquilem} For every $\psi\in H^\infty$, the Toeplitz operator $T_{\ov\psi}$ given by 
$$T_{\ov\psi}f:=P_+(\ov\psi f),\qquad f\in\bmoa,$$
maps $\bmoa$ boundedly into itself.
\end{lem}

To prove this, it suffices to observe (in the spirit of \cite{H}) that $T_{\ov\psi}$ is the adjoint of the multiplication operator $g\mapsto\psi g$, which is obviously bounded on $H^1$.

\par Before proceeding, we need to introduce a bit of notation. Namely, with an inner function $\th$ and a number $\eps\in(0,1)$ we associate the sublevel set 
$$\Omte:=\{z\in\D:\,|\th(z)|<\eps\}.$$
The following result is a restricted version of \cite[Theorem 1]{DAJM}.

\begin{lem}\label{lem:secondlem} Suppose that $f\in\bmoa$ and $\th$ is an inner function. Then $f\ov\th\in\bmo$ if and only if 
\begin{equation}\label{eqn:bddlevset}
\sup\{|f(z)|:\,z\in\Omte\}<\infty
\end{equation}
for some (or every) $\eps$ with $0<\eps<1$.
\end{lem}

Next, we recall a remarkable maximum principle for $K^2_\th$ functions that was established by Cohn in \cite{C2}. 

\begin{lem}\label{lem:cohnmax} Let $\th$ be inner, and suppose $f\in K^2_\th$ is a function that satisfies \eqref{eqn:bddlevset} for some $\eps\in(0,1)$. Then $f\in H^\infty$.
\end{lem}

Our last lemma reproduces yet another result of Cohn (see \cite[p.\,737]{C1}), which characterizes the inner functions $\th$ with the property that $K_{*\th}$ contains only bounded functions. This characterization is, in turn, a consequence of Hru\v s\v cev and Vinogradov's earlier work from \cite{HV} on the multipliers of Cauchy type integrals.

\begin{lem}\label{lem:kbmoinf} Let $\th$ be an inner function. Then $K_{*\th}=K^\infty_\th$ if and only if $\th$ is a Frostman Blaschke product.
\end{lem}

We also refer to \cite[Theorem 1.7]{DIUMJ07} for a refinement of this result in terms of $\text{\rm inn}(K_{*\th})$, the set of inner factors for functions from $K_{*\th}$.

\section{Proof of Theorem \ref{thm:intsmooth}}

We shall only give a detailed proof of part (a), the other cases being similar. Since $\mathcal W=\{w_k\}\in\ell^2_1(\mathcal Z)$, we know that there exists a unique $f\in K^2_B$ such that $f\big|_{\mathcal Z}=\mathcal W$. Therefore, in order that
\begin{equation}\label{eqn:inclip}
\mathcal W\in\left(K^2_B\cap\La^\al\right)\big|_{\mathcal Z}
\end{equation}
it is necessary and sufficient that 
\begin{equation}\label{eqn:finla}
f\in\La^\al.
\end{equation}
To find out when the latter condition holds, we apply Lemma \ref{lem:firstlem}, part (a), to the function $g:=\ov z\ov fB$ in place of $f$. (Note that $g\in H^2$ because $f\in K^2_B$.) This tells us that $P_-(\ov Bg)\in\La^\al$ if and only if
\begin{equation}\label{eqn:gzko}
|g(z_k)|=O\left((1-|z_k|)^\al\right),\qquad k\in\N.
\end{equation}
On the other hand, 
$$P_-(\ov Bg)=P_-(\ov z\ov f)=\ov z\ov f,$$
and it is clear that the function $\ov z\ov f$ belongs to $\La^\al$ if and only if $f$ does. Thus, we may rephrase \eqref{eqn:finla} as \eqref{eqn:gzko}. 
\par To arrive at a further---and definitive---restatement of \eqref{eqn:gzko}, we need to express the numbers $g(z_k)$ in terms of $\mathcal W$. For $z\in\D$, Cauchy's formula yields
$$g(z)=\int_{\T}\f{g(\ze)}{1-\ov\ze z}\,dm(\ze).$$
Consequently, 
\begin{equation*}
\begin{aligned}
\ov{g(z_k)}&=\int_{\T}\f{\ov{g(\ze)}}{1-\ze\ov z_k}\,dm(\ze)
=\int_{\T}\f{\ze f(\ze)\ov{B(\ze)}}{1-\ze\ov z_k}\,dm(\ze)\\
&=\f1{2\pi i}\int_{\T}\f{f(\ze)}{B(\ze)\cdot(1-\ze\ov z_k)}\,d\ze.
\end{aligned}
\end{equation*}
Computing the last integral by residues, while recalling that $f(z_j)=w_j$, we find that 
\begin{equation}\label{eqn:bargak}
\ov{g(z_k)}=\sum_j\f{w_j}{B'(z_j)\cdot(1-z_j\ov z_k)}=\widetilde w_k 
\end{equation}
for each $k\in\N$. (To justify the application of the residue theorem, one may begin by evaluating the integral over the circle $r_n\T$, where $\{r_n\}\subset(0,1)$ is a suitable sequence tending to $1$, and then pass to the limit as $n\to\infty$.)
\par Finally, we use \eqref{eqn:bargak} to rewrite \eqref{eqn:gzko} in the form 
\begin{equation}\label{eqn:wtwk}
|\widetilde w_k|=O\left((1-|z_k|)^\al\right),\qquad k\in\N.
\end{equation}
The equivalence of \eqref{eqn:inclip} and \eqref{eqn:wtwk} is thereby established, proving the $\La^\al$ part of the theorem. 
\par The remaining statements (i.e., those involving $\bmo$, $G_\al$ and $\mathcal L^p_s$) are proved similarly, by combining the appropriate parts of Lemma \ref{lem:firstlem} with identity \eqref{eqn:bargak}.

\section{Proofs of Proposition \ref{prop:nonint} and Theorem \ref{thm:nondual}}

We begin by proving Theorem \ref{thm:nondual}. Once this is done, Proposition \ref{prop:nonint} will be derived as a corollary. 

\medskip{\it Proof of Theorem \ref{thm:nondual}.} Given an inner function $\th$ distinct from a finite Blaschke product, we want to find a function $g\in K^2_\th\setminus\bmo$ that induces a bounded linear functional on $K^1_\th$. We shall distinguish two cases.

\smallskip{\it Case 1.} Assume that $\th$ is an infinite Frostman Blaschke product. Its zero sequence, say $\mathcal Z=\{z_j\}$, must then have a limit point on $\T$. Of course, nothing is lost by assuming that $\mathcal Z$ clusters at $1$. Now let 
$$\ph(z):=\log(1-z),$$
where \lq\lq log" stands for the holomorphic branch of the logarithm that lives on the right half-plane and satisfies $\log1=0$. We have $\ph\in\bmoa$ (because $\text{\rm Im}\,\ph\in L^\infty$), so the corresponding linear functional acts boundedly on $H^1$ and hence on $K^1_\th$. Clearly, the same functional on $K^1_\th$ is also induced, in a similar manner, by the function 
$$g:=\th P_-(\ov\th\ph),$$
which is the orthogonal projection (in $H^2$) of $\ph$ onto $K^2_\th$. Precisely speaking, the functional
$$f\mapsto\int_\T f\ov g\,dm\left(=\int_\T f\ov\ph\,dm\right),\qquad f\in K^2_\theta,$$
extends continuously to $K^1_\theta$. 
\par We know that $g\in K^2_\th$, and to conclude that $g$ does the job, we only need to check that 
\begin{equation}\label{eqn:gnotinbmo}
g\notin\bmo.
\end{equation}
To this end, observe first that $\sup_j|\ph(z_j)|=\infty$ and hence, {\it a fortiori},
$$\sup\{|\ph(z)|:\,z\in\Omte\}=\infty$$
for every $\eps\in(0,1)$. By Lemma \ref{lem:secondlem}, this implies that $\ov\th\ph\notin\bmo$. On the other hand, 
$$\ov\th\ph=P_-(\ov\th\ph)+P_+(\ov\th\ph),$$
where the last term, $P_+(\ov\th\ph)$, is in $\bmoa(\subset\bmo)$ thanks to Lemma \ref{lem:sesquilem}. It follows readily that $P_-(\ov\th\ph)\notin\bmo$. In particular, $P_-(\ov\th\ph)\notin L^\infty$ (just note that $L^\infty\subset\bmo$). Equivalently, the function $\th P_-(\ov\th\ph)=g$ is not in $L^\infty$. 
\par Now, if $g$ were in $\bmo$, then we would have $g\in K_{*\th}$; and since our current assumption on $\th$ yields $K_{*\th}=K^\infty_\th$ (in accordance with Lemma \ref{lem:kbmoinf}), $g$ would have to be bounded, which it is not. This proves \eqref{eqn:gnotinbmo}. 

\smallskip{\it Case 2.} Assume that $\th$ is not a Frostman Blaschke product. This time, using Lemma \ref{lem:kbmoinf} again, we can find an unbounded function $h\in K_{*\th}$. We have then $\widetilde h:=\ov z\ov h\th\in K^2_\th$, and we go on to claim that $\widetilde h\notin\bmo$. (Here and below, the \lq\lq tilde operation" \eqref{eqn:antilinisom} is being used repeatedly.) Indeed, if $\widetilde h$ were in $\bmo$, then so would be $h\ov\th$, and Lemma \ref{lem:secondlem} would tell us that 
$$\sup\{|h(z)|:\,z\in\Omte\}<\infty$$
for some (any) $\eps\in(0,1)$. This, however, would imply that $h\in H^\infty$ by virtue of Lemma \ref{lem:cohnmax}, whereas $h$ is actually unbounded by assumption. 
\par Now we know that $\widetilde h\in K^2_\th\setminus\bmo$, and we proceed by showing that $\widetilde h$ generates a continuous linear functional on $K^1_\th$. This will allow us to conclude that $\widetilde h$ is eligible as $g$ (the function we are looking for), and the proof will be complete. 
\par Given $f\in K^2_\th$, we have the elementary identity $\ov f\widetilde h=\widetilde f\ov h$. Recalling also the facts that $\widetilde f\in H^2(\subset H^1)$ and $h\in\bmoa$, we use the duality relation $(H^1)^*=\bmoa$ to infer that
\begin{equation*}
\begin{aligned}
\left|\int_{\T}f\,\ov{\widetilde h}\,dm\right|&=
\left|\int_{\T}\ov f\,\widetilde h\,dm\right|=\left|\int_{\T}\widetilde f\,\ov h\,dm\right|\\
&\le C\|\widetilde f\|_1\|h\|_*=C\|f\|_1\|h\|_*
\end{aligned}
\end{equation*}
with some absolute constant $C>0$. Consequently, for $g=\widetilde h$, the (densely defined) functional
\begin{equation}\label{eqn:fgfunct}
f\mapsto\int_{\T}f\ov g\,dm
\end{equation}
is indeed continuous on $K^1_\th$, and we are done. 
\qed

\bigskip{\it Proof of Proposition \ref{prop:nonint}.} Assuming that $B=B_\mathcal Z$ is an infinite Blaschke product with zeros $\mathcal Z=\{z_j\}$, where the $z_j$'s are pairwise distinct, we want to find a sequence of values $\mathcal W=\{w_j\}$ in $\ell^\infty_{\log}(\mathcal Z)$ that is not the trace of any $K_{*B}$ function on $\mathcal Z$. Consider, for each $j\in\N$, the function 
$$f_j(z):=\f1{1-\ov z_jz}$$
and note that $f_j\in K^2_B$. Observe also that
\begin{equation}\label{eqn:fjlog}
\|f_j\|_1\le M\log\f2{1-|z_j|},\qquad j\in\N,
\end{equation}
for some fixed constant $M>0$. 
\par Now, Theorem \ref{thm:nondual} provides us with a function $g\in K^2_B\setminus\bmo$ such that the associated functional \eqref{eqn:fgfunct}, defined initially for $f\in K^2_B$, acts boundedly on $K^1_B$, say with norm $N_g$. When applied to $f=f_j$, this functional takes the value $\ov{g(z_j)}$; indeed, Cauchy's formula gives
$$\int_\T\ov f_jg\,dm=g(z_j)$$
for each $j$. In conjunction with \eqref{eqn:fjlog}, this yields
\begin{equation}\label{eqn:gzjlog}
|g(z_j)|\le N_g\|f_j\|_1\le MN_g\log\f2{1-|z_j|},\qquad j\in\N.
\end{equation}
Finally, we put 
$$w_j:=g(z_j),\qquad j\in\N.$$
The sequence $\mathcal W=\{w_j\}$ is then in $\ell^\infty_{\log}(\mathcal Z)$, as \eqref{eqn:gzjlog} shows, while
\begin{equation}\label{eqn:lasteq}
\mathcal W\notin K_{*B}\big|_\mathcal Z
\end{equation}
as required. To verify \eqref{eqn:lasteq}, it suffices to note that $g$ is the {\it only} function in $K^2_B$ that interpolates $\mathcal W$ on $\mathcal Z$ (indeed, a $K^2_B$ function is uniquely determined by its trace on $\mathcal Z=B^{-1}(0)$), whereas $g\notin\bmo$. The proof is complete. 
\qed

\medskip{\it Remark.} We have seen above that if $g\in K^2_B$, with $B=B_{\mathcal Z}$, and if the functional \eqref{eqn:fgfunct} is continuous on $K^1_B$, then $g\big|_{\mathcal Z}\in\ell^\infty_{\log}(\mathcal Z)$. Now, if $\mathcal Z$ has the $\bmoa$-interpolating property \eqref{eqn:bmoazlog}, then {\it every} sequence $\mathcal W$ in $\ell^\infty_{\log}(\mathcal Z)$ is actually writable as $g\big|_{\mathcal Z}$ for some $g\in K^2_B$ that induces a continuous linear functional on $K^1_B$. (To see why, take $G\in\bmoa$ with $G\big|_{\mathcal Z}=\mathcal W$ and then put $g=BP_-(\ov BG)$, so that $g$ is the orthogonal projection of $G$ onto $K^2_B$.) Of course, things become different if condition \eqref{eqn:bmoazlog} is dropped. For instance, there are interpolating sequences $\mathcal Z$ for which $\ell^\infty_{\log}(\mathcal Z)\not\subset\ell^2_1(\mathcal Z)$; and if this is the case, then no sequence in $\ell^\infty_{\log}(\mathcal Z)\setminus\ell^2_1(\mathcal Z)$ is the trace of any $H^2$ function on $\mathcal Z$.

\end{document}